\documentclass[12pt]{amsart}
\usepackage{amsmath,amscd,amssymb,amsfonts}
\newcommand{\h}{\hbox}
\newcommand{\q}{\quad}

\newcommand{\nin}{\par\noindent}
\newcommand{\bs}{\par\bigskip}
\newcommand{\ms}{\par\medskip}
\newcommand{\sk}{\par\smallskip}
\newcommand{\mopl}{\h{$\bigoplus$}}
\newcommand{\msum}{\h{$\sum$}}
\newcommand{\al}{\alpha}
\newcommand{\C}{{\mathbf C}}
\newcommand{\D}{{\mathbf D}}
\newcommand{\De}{\Delta}
\newcommand{\DS}{\Delta^*}
\newcommand{\dd}{\partial}
\newcommand{\f}{{}\,\overline{\!f}}

\newcommand{\F}{{\mathcal F}}
\newcommand{\FD}{{\mathcal F}_{\Delta^*}}
\newcommand{\Ft}{\widetilde{\mathcal F}}
\newcommand{\FDt}{\widetilde{{\mathcal F}|_{\Delta^*}}}
\newcommand{\ga}{\gamma}
\newcommand{\HM}{H_{\mathcal M}}
\newcommand{\Hs}{H_s}
\newcommand{\Hsc}{H_s^c}
\newcommand{\Hi}{H_{\infty}}
\newcommand{\Hic}{H_{\infty}^c}
\newcommand{\Hiu}{H_{\infty,1}}
\newcommand{\Hiuc}{H_{\infty,1}^c}
\newcommand{\Hinu}{H_{\infty,\ne 1}}
\newcommand{\Hinuc}{H_{\infty,\ne1}^c}
\newcommand{\Htu}{\widetilde{H}_1}
\newcommand{\Hu}{H_1}
\newcommand{\Huq}{H_{1,\Q}}
\newcommand{\iy}{\infty}
\newcommand{\M}{{\mathcal M}}
\newcommand{\MM}{\widetilde{\mathcal M}}
\newcommand{\PP}{{\mathbf P}}
\newcommand{\Q}{{\mathbf Q}}
\newcommand{\R}{{\mathbf R}}
\newcommand{\So}{{}\,\overline{\!S}}
\newcommand{\V}{{}_0V}
\newcommand{\X}{\overline{X}}
\newcommand{\Y}{\overline{Y}}
\newcommand{\Yt}{\widetilde{Y}}
\newcommand{\Z}{{\mathbf Z}}
\newcommand{\Gr}{{\rm Gr}}
\newcommand{\IC}{{\rm IC}}
\newcommand{\Ext}{{\rm Ext}}
\newcommand{\Hom}{{\rm Hom}}
\newcommand{\Ker}{{\rm Ker}}
\newcommand{\MHM}{{\rm MHM}}
\newcommand{\MHS}{{\rm MHS}}
\newcommand{\Coker}{{\rm Coker}}
\newcommand{\Perv}{{\rm Perv}}
\newcommand{\supp}{{\rm supp}}
\newcommand{\into}{\hookrightarrow}
\newcommand{\bl}{\bigl}
\newcommand{\br}{\bigr}
\newcommand{\ssb}{\raise.15ex\h{${\scriptscriptstyle\bullet}$}}
\newcommand{\ssc}{\,\raise.15ex\h{${\scriptstyle\circ}$}\,}
\newcommand{\simto}{\buildrel{\sim}\over\longrightarrow}

\begin{document}
\title[Weight filtration of limit mixed Hodge structure]
{Weight filtration of the limit mixed Hodge structure at infinity
for tame polynomials}
\author{Alexandru Dimca}
\address{Institut Universitaire de France et
Laboratoire J.A.\ Dieudonn\'e, UMR du CNRS 7351,
Universit\'e de Nice-Sophia Antipolis, Parc Valrose,
06108 Nice Cedex 02, France}
\email{Alexandru.DIMCA@unice.fr}
\author{Morihiko Saito}
\address{RIMS Kyoto University, Kyoto 606-8502 Japan}
\email{msaito@kurims.kyoto-u.ac.jp}
\begin{abstract}
We give three new proofs of a theorem of C. Sabbah asserting that the
weight filtration of the limit mixed Hodge structure at infinity of
cohomologically tame polynomials coincides with the monodromy filtration
up to a certain shift depending on the unipotent or non-unipotent
monodromy part.
\end{abstract}
\maketitle

\centerline{\bf Introduction}

\bs\nin
Let $X=\C^n\,\,(n\ge 2)$, and $S=\C$.
Let $f:X\to S$ be a cohomologically tame polynomial map in the sense of
[Sab3], i.e.\ there is a relative compactification $\f:\X\to S$ of $f$
such that $\f$ is proper and the support of $\varphi_{\f-c}\R j_*\Q_X$
is contained in the affine space $X$ (hence discrete) for any $c\in\C$,
where $j:X\into\X$ denotes the inclusion.
\sk
Set $m=n-1$, and $X_s:=f^{-1}(s)\subset X$ for $s\in S$.
There is a non-empty Zariski-open subset $S'$ of $S$ such that $X_s$
is smooth for $s\in S$ and the $H^i(X_s,\Q)$ form a local system on
$S'$. For $s\in S'$, we have moreover
$$H^i(X_s,\Q)=0\,\,\,\h{if}\,\,\,i\ne 0,\,m,\q
H_c^i(X_s,\Q)=0\,\,\,\h{if}\,\,\,i\ne m,\,2m,$$
and $H^0(X_s,\Q)=\Q$, $H_c^{2m}(X_s,\Q)=\Q(-m)$.
These follow from the discreteness of the support of
$\varphi_{\f-c}\R j_*\Q_X$ together with $X=\C^n$ by using the Leray
spectral sequence as in Remark~(1.2) below. For $s\in S'$, set
$$\Hs:=H^m(X_s,\Q),\q \Hsc:=H_c^m(X_s,\Q).$$
By [De1], these have the canonical mixed Hodge structures which are
dual of each other (up to a Tate twist). By the definition of the
weight filtration $W$ in loc.~cit., Th.~3.2.5\,(iii), we have
$\Gr^W_{m+k}\Hs=\Gr^W_{m-k}\Hs=0$ for $k<0$.
In [DS], Th.~0.3, the following was shown:
$$\{\Gr^W_{m+k}\Hs\}_{s\in S'},\,\,
\{\Gr^W_{m-k}\Hsc\}_{s\in S'}\,\,\,\h{are
{\it constant} on $S'$ if $k>0$.}
\leqno(0.1)$$
In fact, the argument in loc.~cit.\ implies that
$\{\Hs/W_m\Hs\}_{s\in S'}$ and $\{W_m\Hsc\}_{s\in S'}$ are constant.
This is closely related to (1.1.3) below.
\sk
Let $\Hi$ be the limit mixed Hodge structure of $\Hs$ for
$s\to\iy$, and similarly for $\Hic$, see [St1], [St2], [StZ].
Set $N:=(2\pi i)^{-1}\log T_u$ with $T_u$ the unipotent part of the
monodromy at infinity. This is an endomorphism of type $(-1,-1)$
of $\Hi$, $\Hic$.
Let $L$ be the filtration on $\Hi$, $\Hic$ induced by the weight
filtration $W$ respectively on $\Hs$, $\Hsc$ for $s\in S'$.
Then the weight filtration $W$ on $\Hi$, $\Hic$ coincides with the
relative monodromy filtration of $(L,N)$, see [StZ].
In particular, $W$ on $\Gr^L_m\Hi$, $\Gr^L_m\Hic$ coincides
with the monodromy filtration shifted by $m$ (i.e.\ with center $m$).
\sk
Let $\Hiu$, $\Hiuc$ respectively denote the unipotent monodromy
part of $\Hi$, $\Hic$,
and similarly for the non-unipotent part $\Hinu$, $\Hinuc$.
By (0.1) we have
$$\Hinu=\Gr^L_m\Hinu,\q\Hinuc=\Gr^L_m\Hinuc.
\leqno(0.2)$$
We thus get the following well-known assertion (see also Appendix of
[MT]):
\ms\nin
{\bf Proposition~1.} {\it With the above notation and assumption,
the weight filtration $W$ on $\Hinu$, $\Hinuc$ coincides with the
monodromy filtration shifted by $m$.}
\ms
In this paper we give three new proofs of the following.
\ms\nin
{\bf Theorem~1} (C.~Sabbah).
{\it With the above notation and assumption, the weight filtration $W$
on $\Hiu$, $\Hiuc$ coincides with the monodromy filtration shifted by
$m+1$ and $m-1$ respectively, and we have the isomorphisms of mixed
Hodge structures for $k\ge 1:$
$$\Gr^W_{m+k}\Hs\cong\Gr^W_{m+k}\Coker(N|\Hiu),\q
\Gr^W_{m-k}\Hsc\cong\Gr^W_{m-k}\Ker(N|\Hiuc),$$
where $\Coker(N|\Hiu)$ is a quotient of $\Hiu(-1)$ and
$\Ker(N|\Hiuc)\subset\Hiuc$.}
\ms
Note that the assertions on $\Hiu$ and $\Hiuc$ are dual of each
other.
The last assertion of Theorem~1 means that the primitive part of the
graded pieces of the monodromy filtration on $\Hiu$ is given by
$\Gr^W_{m+k}\Hs$ for $k\ge 1$, and the coprimitive part for $\Hiuc$ by
$\Gr^W_{m-k}\Hsc$.
\sk
Theorem~1 was first obtained by C.~Sabbah as a corollary of
[Sab3], Th.~13.1 where he uses a theory of Fourier transformation,
Brieskorn lattices, and spectra at infinity, which was developed by
him (see also [Sab1], [Sab2]).
Recently another proof has been given also by him in Appendix of [MT]
without using Brieskorn lattices or spectra at infinity, but using
Fourier transformation where irregular $D$-modules inevitably appear.
\sk
It does not seem, however, that the above theory is absolutely
indispensable for the proof of Theorem~1.
In fact, the theorem was almost proved in [DS] where the following was
shown (see also [Di1], 4.3--5):
\ms\nin
{\bf Theorem~2} ([DS], Th.~0.3). {\it With the above notation and
assumption,
let $\nu_k$ and $\nu'_k$ denote the number of Jordan blocks of size
$k$ for the monodromy on $\Hiu$ and $\Gr^L_m\Hiu$ respectively.
Let $s\in S'$. Then}
$$\nu_k=\dim\Gr^W_{m+k}\Hs,\,\,\,\,\nu'_k=\nu_{k+1}
\q\h{for any}\,\,\,k\ge 1.$$
\ms
We give the first proof of Theorem~1 in this paper by showing that
Theorem~2 implies Theorem~1 using some lemma of linear algebra,
see (1.3--4) below.
\sk
The second proof of Theorem~1 in this paper uses a geometric argument
together with duality, and is quite different from (and perhaps more
intuitive than) the one in the proof of Th.~0.3 in [DS].
It is finally reduced to the following:
\ms\nin
{\bf Proposition~2.} {\it Let $\iota_s:\Hsc\to\Hs$ be the
natural morphism of mixed Hodge structures for $s\in S'$.
Then it induces an isomorphism of Hodge structures}
$$\Gr^W_m\iota_s:\Gr^W_m\Hsc\simto\Gr^W_m\Hs\q\h{for}\,\,\,s\in S'.$$
\ms
We give two proofs of Proposition~2 in this paper.
One proof uses semisimplicity of pure Hodge modules together with
a certain property of the mixed Hodge module $H^0f_*(\Q_{h,X}[n])$
coming from the cohomologically tame condition.
(For $\Q_{h,X}$, see (1.1) below.)
Another proof uses positivity of the polarization on the primitive
cohomology of a compact K\"ahler manifold together with Hironaka's
resolution of singularities.
\sk
The third proof of Theorem~1 in this paper is given as a corollary of
Theorem~3 below, which holds for any pure Hodge module $\M$ of weight
$n$ on $S$ without a constant direct factor, and was proved by
C.~Sabbah in Appendix of [MT]. For a bounded complex of mixed Hodge
modules $\M^{\ssb}$ on a complex algebraic variety $X$ in general, we
denote the mixed Hodge structure $H^j(a_X)_*\M^{\ssb}$ by
$H^j(X,\M^{\ssb})$ (using a remark before (1.1.5) below), where
$a_X:X\to pt$ is the canonical morphism. This notation is compatible
with that of the cohomology of the underlying $\Q$-complex.
\ms\nin
{\bf Theorem~3} (C.~Sabbah).
{\it Let $\M$ be a pure Hodge module of weight $n$
on $S$ having no constant direct factor, i.e.\ $H^{-1}(S,\M)=0$.
Set $H_1:=\psi_{1/t,1}\M$.
Then the weight filtration on $H_1$ coincides with the monodromy
filtration shifted by $n$, and the $N$-primitive part $P\Gr^W_kH_1$
is given by $\Gr^W_kH^0(S,\M)$ for any $k$ where they vanish
unless $k\ge n$.}
\ms
The proof of Theorem~3 in this paper uses the notion of representative
functor and the universal extension $\MM$ of a pure Hodge module $\M$
by a constant mixed Hodge module (see (3.1) below), but Fourier
transformation is not used. Note that $\MM$ was defined in loc.~cit.\
by using a sheaf-theoretic operation explicitly.
By the property of the universal extension we have
$$\MM/\M=a_X^*H^0(S,\M)[1],\q
H^{-1}(S,\MM/\M)=H^0(S,\M).$$
The proof of Theorem~3 is reduced to the comparison between the global
universal extension on $S$ and the local one on a neighborhood of
$\infty\in\PP^1$ for the underlying perverse sheaves, where some
argument is similar to the one in the proof of [DS], Th.~0.3.
\sk
We thank the referee for helping us to improve the paper.
\sk
The first named author was partially supported by the grant
ANR-08-BLAN-0317-02 (SEDIGA).
The second named author is partially supported by Kakenhi 21540037.
\sk
In Section~1 we explain some basics on cohomologically tame
polynomials, and give the proof of Theorem~1 using Theorem~2
after showing Lemma~(1.3).
In Section~2 we give two proofs of Proposition~2, and then a
geometric proof of Theorem~1 after showing Lemma~(2.3).
In Section~3 we explain the notion of a universal extension by a
constant sheaf, and then prove Theorem~3 which implies Theorem~1. 

\bs\bs
\centerline{\bf 1. Cohomologically tame polynomials}
\bs\nin
In this section we explain some basics on cohomologically tame
polynomials, and give the proof of Theorem~1 using Theorem~2
after showing Lemma~(1.3).
\ms\nin
{\bf 1.1.~Some basics on cohomologically tame polynomials.}
Set $X=\C^n\,\,(n\ge 2)$, and $S=\C$. Let $f:X\to S$ be a polynomial
map, and $\f:\X\to S$ be an algebraic compactification of $f$
(i.e. $\f$ is proper). Let $j:X\into\X$ be the inclusion.
Note that $\R j_*\Q_X[n]$ is a perverse sheaf since $j$ is an
affine open immersion, see [BBD].
The intersection complex $\IC_{\X}\Q$ is a subobject of the
perverse sheaf $\R j_*\Q_X[n]$ (see loc.~cit.) and the vanishing
cycle functor $\varphi_{\f-c}$ (see [De2]) is an exact functor of
perverse sheaves (up to a shift).
So we get the first inclusion of
$$\supp\,\varphi_{\f-c}\IC_{\X}\Q\subset
\supp\,\varphi_{\f-c}\R j_*\Q_X[n]=
\supp\,\varphi_{\f-c}\R j_!\Q_X[n].
\leqno(1.1.1)$$
For the last isomorphism, we have the relation
$\D\ssc\R j_*=\R j_!\ssc\D$ and the compatibility of $\varphi_{\f-c}$
with the dualizing functor $\D$, i.e.\ $\D\ssc\varphi=\varphi\ssc\D$
where a Tate twist may appear depending on the eigenvalue of the
monodromy, see e.g. [Sai1], 5.2.3.
\sk
Assume now that $\supp\,\varphi_{\f-c}\R j_*\Q_X$ is contained in the
affine space $X$ (hence discrete) for any $c\in\C$. This means that
$f$ is a {\it cohomologically tame} polynomial in the sense of [Sab3].
Consider the canonical morphisms
$$\R j_!\Q_X[n]\to\R j_*\Q_X[n],\q
\R j_!\Q_X[n]\to\IC_{\X}\Q,\q
\IC_{\X}\Q\to\R j_*\Q_X[n].
\leqno(1.1.2)$$
Since the restrictions of these morphisms to $X$ are the identity
morphisms, and $\varphi$ commutes with the direct images by proper
morphisms, we get the following.
\begin{itemize}
\item[(1.1.3)]
The direct images by $\f$ of the mapping cones of the canonical
morphisms in (1.1.2) are isomorphic to direct sums of shifted constant
sheaves on $S$.
\end{itemize}
Indeed, the above properties imply constancy of the cohomology
sheaves. Then (1.1.3) follows from the vanishing of
$\Ext^i(\Q_S,\Q_S)$ for $i>0$ by using the canonical filtration
$\tau_{\le k}$ in [De1].
\ms
Let $\Q_{h,X}$ denote the object in $D^b\MHM(X)$
(the bounded derived category of mixed Hodge modules on $X$)
which is uniquely characterized by the following two conditions:
Its underlying $\Q$-complex is $\Q_X$, and its $0$-th cohomology
$H^0(X,\Q_{h,X}):=H^0(a_X)_*\Q_{j,X}$ has weight $0$,
see [Sai2], 4.4.2. We have
\begin{itemize}
\item[(1.1.4)]
Replacing $\Q_X$ with $\Q_{h,X}$ in (1.1.2), the assertion (1.1.3) holds
in $D^b\MHM(S)$.
\end{itemize}
Indeed, any admissible variation of mixed Hodge structure $\M$ on $S$
is a constant variation, see e.g.\ [StZ], Prop.~4.19. (This follows from
the existence of the canonical mixed Hodge structure on $H^0(S,\M)$ by
using the restriction morphisms $H^0(S,\M)\to\M_s$ which are morphisms
of mixed Hodge structures for $s\in S$, where $\M_s$ is the pull-back
of $\M$ by the inclusion $\{s\}\into S$, and $\M[1]$ is a mixed Hodge
module on $S$.) So we get $\M=a_S^*H$ for $H\in\MHS$ (the category of
graded-polarizable mixed Hodge $\Q$-structures in [De1]), where
$a_S:S\to pt$ is the canonical morphism. Note that $\MHM(pt)$ is
naturally identified with $\MHS$, see [Sai2]. We have moreover for $i>1$
$$\Ext^i_{\MHM(S)}(a_S^*H,a_S^*H')=\Ext^i_{\MHS}(H,(a_S)_*a_S^*H')=
\Ext^i_{\MHS}(H,H')=0,
\leqno(1.1.5)$$
since $\Ext^i=0\,\,(i>1)$ in $\MHS$ by a well-known corollary of a
theorem of Carlson [Ca] (which implies the right-exactness of the
functor $\Ext^1_{\MHS}(\Q,*)$). So (1.1.4) follows by using the
canonical filtration $\tau_{\le k}$ in [De1].

\ms\nin
{\bf 1.2.~Remark.} Let $f:X\to S$ be as in the beginning of (1.1).
We have the Leray spectral sequence in $\MHS$:
$$E_2^{i,j}=H^i(S,H^jf_*(\Q_{h,X}[n]))\Rightarrow H^{i+j+n}(X,\Q),
\leqno(1.2.1)$$
using the canonical filtration $\tau_{\le k}$ as in [De1].
(This will be used later.)
\sk
We have $E_2^{i,j}=0$ for $i\notin[-1,0]$ since $S=\C$.
So (1.2.1) degenerates at $E_2$, and we get
$$H^i(S,H^jf_*(\Q_{h,X}[n]))=0\q\h{for}\,\,\,(i,j)\ne(-1,1-n),
\leqno(1.2.2)$$
since $X=\C^n$ and $H^jf_*(\Q_{h,X}[n])=0$ for $j\le -n$
(using the classical $t$-structure).
\sk
Assume $f$ is cohomologically tame (using an appropriate
compactification of $f$ as in (1.1)). Then
$$H^jf_*(\Q_{h,X}[n])\,\,\,\h{is constant for}\,\,\,j\ne 0,
\leqno(1.2.3)$$
by using the exactness of $\varphi$ (up to a shift) together with
the commutativity of $\varphi$ and the direct image under proper
morphisms as in (1.1). So (1.2.2) implies
$$H^jf_*(\Q_{h,X}[n])=0\q\h{unless $\,j=1-n\,$ or $\,0$,}
\leqno(1.2.4)$$

\ms\nin
{\bf 1.3.~Lemma.} {\it Let $V_{\ssb}$ be a finite dimensional graded
$\Q$-vector space with an action of $N$ of degree $-2$, i.e.
$N(V_k)\subset V_{k-2}$. Let $V'_{\ssb}$ be a graded vector subspace
stable by $N$.
Set $V''_{\ssb}:=V_{\ssb}/V'_{\ssb}$. Let $m$ be an integer.
Assume the action of $N$ on $V''_{\ssb}$ vanishes, and
$$N^k:V'_{m+k}\simto V'_{m-k}\q\h{for any}\,\,\,k\ge 1.$$
Set $C'_k:=\Coker(N:V'_{k+2}\to V'_k)$ so that $N$ induces
$\delta_k:V''_{m+k+2}\to C'_{m+k}$.
Let $\nu_k$ be the number of Jordan blocks of size $k$ for the
action of $N$ on $V_{\ssb}$. Then}
$$\nu_{k+1}=\begin{cases}\dim\Coker\,\delta_0+\msum_j\dim\Ker\,
\delta_j&\h{if}\,\,\,k=0,\\
\dim\Coker\,\delta_k+\dim{\rm Im}\,\delta_{k-1}\raise12pt\h{ }
&\h{if}\,\,\,k\ge 1.\end{cases}$$

\ms\nin
{\it Proof.}
Let $\V'_{m+k}$ be the primitive part defined by
$\Ker\,N^{k+1}\subset V'_{m+k}$ for $k\ge 0$.
We have the primitive decomposition
$$V'_{\ssb}=\mopl_{k\ge 0}\bigl(\mopl_{j=1}^kN^j\V'_{m+k}\bigr)
\q\h{with}\q\V'_{m+k}\simto C'_{m+k}.
\leqno(1.3.1)$$
Set
$$n_k=\dim{\rm Im}\,\delta_k.$$
For each $k\ge 0$, there are bases $\{v'_{k,j}\}_j$ of
$\V'_{m+k}\,(=C'_{m+k})$ and $\{v''_{k,j}\}_j$ of $V''_{m+k+2}$
together with lifts $v_{k,j}$ of $v''_{k,j}$ in $V_{m+k+2}$ such that
$$Nv_{k,j}=\begin{cases}v'_{k,j}&\h{if}\,\,\,1\le j\le n_k,\\
0&\h{otherwise}.\end{cases}
\leqno(1.3.2)$$
Indeed, by the definition of $\delta_k$, the assertion is trivial
if we consider the equality modulo $NV'_{\ssb}$,
i.e.\ if we add the term $+Nu_{k,j}$ for some $u_{k,j}\in V'_{m+k+2}$
on the right-hand side of (1.3.2).
Then we can replace the lift $v_{k,j}$ of $v''_{k,j}$ with
$v_{k,j}-u_{k,j}$, and (1.3.2) is proved.

The assertion of Lemma~(1.3) now follows from (1.3.1) and (1.3.2).

\ms\nin
{\bf 1.4.~Remark.}
Assume the morphisms $\delta_k:V''_{m+k+2}\to C'_{m+k}$ in Lemma~(1.3)
are bijective for any $k\ge 0$. Then (1.3.2) in the proof of
Lemma~(1.3) implies that the the primitive decomposition of
$V'_{m+\ssb}$ can be lifted to that of $V_{m+1+\ssb}$.

\ms\nin
{\bf 1.5.~Proof of Theorem~1 using Theorem~2.}
We show the assertion for $\Hiu$ since this implies the assertion
for $\Hiuc$ by duality.
We can replace $\Hiu$ with the graded pieces $\Gr^W_{\ssb}\Hiu$
in order to define $\nu_k$, $\nu'_k$, since $W$ is strictly
compatible with $N^k$ for any $k\ge 1$.
We then apply Lemma~(1.3) to
$$V_k=\Gr^W_k\Hiu,\q V'_k=\Gr^W_k\Gr^L_m\Hiu,\q
V''_k=\begin{cases}\Gr^L_k\Hiu&\h{if}\,\,\,k>m,\\
0&\h{if}\,\,\,k\le m.\end{cases}$$
Here $\Gr^W_j\Gr^L_k\Hiu=0$ for $j\ne k$ if $k>m$, since
$N=0$ on $\Gr^L_k\Hiu$ for $k>m$ by (0.1).
\sk
Using the primitive decomposition (1.3.1), we get
$$\nu'_{k+1}=\dim C'_{m+k}\q\h{for}\,\,\,k\ge 0.
\leqno(1.5.1)$$
We show that Theorem~2 together with Lemma~(1.3) imply the isomorphism
$$\delta_k:V''_{m+k+2}\simto C'_{m+k}\q\h{for}\,\,\,k\ge 0.
\leqno(1.5.2)$$
By (1.5.1) the surjectivity of $\delta_k$ is equivalent to
$$\nu'_{k+1}=\dim{\rm Im}\,\delta_k\q\h{for}\,\,\,k\ge 0,
\leqno(1.5.3)$$
and we have by Theorem~2 and Lemma~(1.3)
$$\nu'_{k+1}=\nu_{k+2}=\dim\Coker\,\delta_{k+1}+
\dim{\rm Im}\,\delta_k\q\h{for}\,\,\,k\ge 0.$$
So (1.5.3) follows by decreasing induction on $k\ge 0$.

As for the injectivity of $\delta_k$, we get by Lemma~(1.3) together
with the surjectivity of $\delta_0$
$$\nu_1=\dim V''_{m+1}+\msum_{k\ge 0}\dim\Ker\,\delta_k,$$
since $\delta_{-1}$ vanishes.
We have moreover $\nu_1=\dim V''_{m+1}$ by Theorem~2.
So the injectivity of $\delta_k\,\,(k\ge 0)$ follows.
Thus (1.5.2) is proved.
\sk
Then the primitive decomposition of $V'_{m+\ssb}$ can be lifted to
that of $V_{m+1+\ssb}$ as is noted in Remark~(1.4).
We have the last assertion of Theorem~1 since the
$\delta_k\,\,(k\ge 0)$ underlie isomorphisms of mixed Hodge structures.
We thus get the first proof of Theorem~1 in this paper.

\bs\bs
\centerline{\bf 2. Geometric proof of Theorem~1}
\bs\nin
In this section we give two proofs of Proposition~2, and then a
geometric proof of Theorem~1 after showing Lemma~(2.3).

\ms\nin
{\bf 2.1.~One proof of Proposition~2.}
Consider the following morphisms of mixed Hodge modules on $S$:
$$\M_!:=H^0\f_*(j_!\Q_{h,X}[n])\buildrel{u'}\over\to
H^0\f_*\IC_{\X}\Q_h\buildrel{v'}\over\to
\M_*:=H^0\f_*(j_*\Q_{h,X}[n]).
\leqno(2.1.1)$$
These are induced by the canonical morphisms of mixed Hodge modules
on $\X$ whose underlying morphisms are as in (1.1.2):
$$j_!\Q_{h,X}[n]\buildrel{u}\over\to
\IC_{\X}\Q_h\buildrel{v}\over\to j_*\Q_{h,X}[n].
\leqno(2.1.2)$$
By (1.1.4) the kernel and cokernel of $u'$ and $v'$ are constant mixed
Hodge modules on $S$.
\sk
By the formalism of mixed Hodge modules (see e.g. [Sai2], 2.26) we have
$$\Gr^W_{n+k}(j_!\Q_{h,X}[n])=\Gr^W_{n-k}(j_*\Q_{h,X}[n])=0
\q\h{if}\,\,\,k>0,
\leqno(2.1.3)$$
and moreover
$$\Gr^W_n(j_!\Q_{h,X}[n])=\Gr^W_n(j_*\Q_{h,X}[n])=\IC_{\X}\Q_h.
\leqno(2.1.4)$$
(Indeed, for the last assertion, we use
${\rm Hom}(\M',j_*\Q_{h,X})={\rm Hom}(j^*\M',\Q_{h,X})=0$ for any
mixed Hodge module $\M'$ supported on $\X\setminus X$, and similarly
for the dual assertion.)
Note that the weight filtration $W$ on $\M_!$ and $\M_*$ is induced
by the weight filtration $W$ on $j_!\Q_{h,X}[n]$ and $j_*\Q_{h,X}[n]$
respectively via the weight spectral sequence, see [Sai2], Prop.~2.15.
Moreover, this $W$ induces the weight filtration $W$ on $H_s$, $H_s^c$
if $s$ is in a sufficiently small non-empty Zariski-open subset of $S$
(since $W$ is independent of the choice of a compactification).
\sk
By (1.1.4) we have constancy of the kernel and the cokernel of the
canonical morphism
$$\Gr^W_n\M_!\to \Gr^W_n\M_*.
\leqno(2.1.5)$$
Assume the cokernel is nonzero. We have $\Gr^W_k\M_*=0$ for $k<n$ by
(2.1.3). So there is a nontrivial constant Hodge submodule in $\M_*$
by semisimplicity of pure Hodge modules applied to $\Gr^W_n\M_*$.
(The latter property follows from polarizability of pure Hodge
modules by using [Sai2], Th.~3.21.)
However, this contradicts the property that $H^{-1}(S,\M_*)=0$ which
follows from the condition that $X=\C^n$ by using the Leray spectral
sequence as in Remark~(1.2). So we get the surjectivity.
For the injectivity we apply the dual argument.
Restricting over a sufficiently general $s\in S$, we then get the desired
isomorphism.

\ms\nin
{\bf 2.2.~Another proof of Proposition~2.} Set $Y=X_s$. More generally,
let $Y$ be a smooth variety which is the complement of an ample
effective divisor $E$ on a projective variety $\Y$. Under this
assumption, we show the bijectivity of the canonical morphism
$$Gr^W_mH_c^m(Y)\to\Gr^W_mH^m(Y).
\leqno(2.2.1)$$
We have a smooth projective compactification $\Yt$ of $Y$ such that
$D:=\Yt\setminus Y$ is a divisor with simple normal crossings.
This is obtained by using Hironaka's resolution
$\sigma:(\Yt,D)\to(\Y,E)$ which is a projective morphism.
Let $D_{\sigma}$ be a relatively ample divisor for $\sigma$. We may
replace $D_{\sigma}$ with $D_{\sigma}-\sigma^*\sigma_*D_{\sigma}$
so that its support is contained in $D$.
Then $k\sigma^*E+D_{\sigma}$ is an ample divisor on $\Yt$ for $k\gg 0$.
Since its support is contained in $D$, it is a linear combination of
the irreducible components $D_i$ of $D$.
\sk
By Deligne's construction of mixed Hodge structure on $H^i(Y)$
(see [De1]) together with duality, we have
$$\aligned Gr^W_mH_c^m(Y)&=\Ker\bl(H^m(\Yt)\to\mopl_i\,H^m(D_i)\br),\\
Gr^W_mH^m(Y)&=\Coker\bl(\mopl_i\,H^{m-2}(D_i)(-1)\to H^m(\Yt)\br).
\endaligned
\leqno(2.2.2)$$
So the assertion is equivalent to the non-degeneracy of the
restriction of the natural pairing on the middle cohomology $H^m(\Yt)$
to the kernel of the morphism $H^m(\Yt)\to\mopl_i\,H^m(D_i)$.
By Hodge theory, it is enough to show that this kernel is contained
in the primitive part with respect to the above ample divisor.
But it is clear since the action of the cohomology class of each $D_i$
is given by composing the restriction morphism
$H^{\ssb}(\Yt)\to H^{\ssb}(D_i)$ with its dual.
So the assertion follows.
This finishes another proof of Proposition~2.

\ms
For the geometric proof of Theorem~1 in this section, we also need the
following.

\ms\nin
{\bf Lemma~2.3.} {\it Let $H$ be a mixed Hodge structure, and $L$ an
increasing filtration on $H$.
Let $N$ be a nilpotent endomorphism of type $(-1,-1)$ of $H$
preserving the filtration $L$.
Assume the relative monodromy filtration $W$ for $(L,N)$ exists, and
$W$ coincides with the weight filtration of the mixed Hodge structure
$H$. Let $m$ be an integer such that $H=L_mH$.
Assume the action of $N$ on $L_{m-1}H$ vanishes, and
$$\aligned\dim\Gr^W_{m-k}(\Ker\,N)=
\dim\Gr^W_{m+k}(\Coker\,N)\q(k\ge 1),\\
\Gr^W_{m-k}(\Ker\,N)=\Gr^W_{m+k}(\Coker\,N)=0\q(k\le 0),\endaligned
\leqno(2.3.1)$$
where $\Coker\,N$ is a quotient of $H(-1)$, and $\Ker\,N\subset H$.
Then $W$ coincides with the monodromy filtration with center $m-1$.}

\bs\nin
{\it Proof.} Set $H':=L_{m-1}H$, $H'':=\Gr^L_mH$.
The hypothesis on the action of $N$ on $H'$ implies that
$$\Gr^W_k\Gr^L_iH'=0\,\,\,(k\ne i),\,\,\,\h{i.e.}\,\,\,
W=L\,\,\,\h{on}\,\,\,H'.
\leqno(2.3.2)$$
Set $H_k:=\Gr^W_kH$, and similarly for $H'_k$, $H''_k$. Set
$$K_k:=\Ker\bl(N:H_k\to H_{k-2}(-1)\br),\q
C_k:=\Coker\bl(N:H_k\to H_{k-2}(-1)\br),$$
and similarly for $K''_k$, $C''_k$.
Note that $\Gr^W_k$ commutes with Ker and Coker by the strict
compatibility of the weight filtration $W$.
\sk
Applying the snake lemma to the action of $N$ on
$0\to H'\to H\to H''\to 0$, we get the following long exact sequence
for any $k\in\Z$
$$0\to H'_k\to K_k\to K''_k\buildrel{\dd}\over\to
H'_{k-2}(-1)\to C_k\to C''_k\to 0.
\leqno(2.3.3)$$
Here $H'_k=K_k=0$ for $k\ge m$, and $C_k=0$ for $k\le m$ by
hypothesis.
\sk
We show the following isomorphisms by decreasing induction on
$k\ge 0$:
$$H'_{m-k}\simto K_{m-k},\q\dd:K''_{m-k}\simto H'_{m-k-2}(-1).
\leqno(2.3.4)$$
Here it is enough to show that $\dim K''_{m-k}=\dim H'_{m-k-2}$,
using the long exact sequence (2.3.3) since the surjectivity of
$\dd$ follows from the vanishing of $C_{k-2}$ for $k\le m$.
\sk
For $k\gg 0$, the assertion trivially holds since all the terms are
zero.
Assume the isomorphisms hold with $k$ replaced by $k+2$.
We have the following equalities for $k\ge 0$:
$$\dim K''_{m-k}=\dim C''_{m+k+2}=\dim C_{m+k+2}=\dim K_{m-k-2}=
\dim H'_{m-k-2}.$$
Indeed, the first equality follows from the property of the
monodromy filtration on $H''$, the second from the long exact sequence
(2.3.3) together with the hypothesis that $H'_{k-2}=0$ for
$k\ge m+2$, the third from the hypothesis (2.3.1) of the lemma,
and the last from the inductive hypothesis.
So the two isomorphisms in (2.3.4) hold for $k\ge 0$.
\sk
We apply Remark~(1.4) to the dual vector space of $H$ and the dual
filtration of $L$ on it. Then (2.3.4) implies that the primitive
decomposition of $\mopl_kH''_k$ with center $m$ can be lifted to the
primitive decomposition of $\mopl_kH_k$ with center $m-1$ under the
surjection $H\to H''$, since $K''_k$ is the coprimitive part of
$H''_k\,(k\le m)$. This finishes the proof of Lemma~(2.3).

\bs\nin
{\bf 2.4.~Proof of Theorem~1.}
We show the assertion for $\Hic$ since that for $\Hi$ follows from
this using duality.
Consider first the following canonical morphisms
$$\R\Gamma_c(S,f_!\Q_{h,X})\buildrel{\al}\over\longrightarrow
\R\Gamma(S,f_!\Q_{h,X})\buildrel{\beta}\over\longrightarrow
\R\Gamma(S,f_*\Q_{h,X}).$$
Set
$$\ga=\beta\ssc\al:\R\Gamma_c(S,f_!\Q_{h,X})\to
\R\Gamma(S,f_*\Q_{h,X}).$$
By the octahedral axiom of the derived category, we get
a distinguished triangle
$$C\bl(\al)\to C(\ga)\to C(\beta)\buildrel{+1}\over\to.
\leqno(2.4.1)$$
By (1.1) the following mapping cone is a direct sum of constant
sheaves on $S$:
$$C(f_!\Q_{h,X}\to f_*\Q_{h,X}\br)=
C\bl(\R\f_*j_!\Q_{h,X}\to\R\f_*j_*\Q_{h,X}\br),$$
and this holds in $D^b\MHM(S)$.
Moreover, the stalk at $s\in S'$ of the mapping cone is given by the
cohomology of the mapping cone
$$C\bl(\R\Gamma_c(X_s,\Q)\to\R\Gamma(X_s,\Q)\br)\q(s\in S'),$$
using the generic base change by the inclusion $\{s\}\into S$.
\sk
We then get the following isomorphisms in the derived category of
graded-polarizable mixed Hodge structures $D^b\MHS$:
$$\aligned C(\al)&=C\bl(N:\psi_{1/t,1}\,f_!\Q_{h,X}\to
\psi_{1/t,1}\,f_!\Q_{h,X}(-1)\br)[-1],\\
C(\beta)&=C\bl(\R\Gamma_c(X_s,\Q)\to\R\Gamma(X_s,\Q)\br)\q(s\in S'),\\
C(\ga)&=\Q\oplus\Q(-n)[3-2n].\endaligned
\leqno(2.4.2)$$
Indeed, the first isomorphism follows from
$$C(j'_!\M\to j'_*\M)=
C\bl(N:\psi_{1/t,1}\M\to\psi_{1/t,1}\M(-1)\br)[-1],$$
for any mixed Hodge module $\M$ on $S$ where
$j':S\into\overline{S}:=\PP^1$ is the inclusion, see [Sai2], 2.24.
(In this paper the nearby and vanishing cycle functors $\psi$,
$\varphi$ for mixed Hodge modules are compatible with those for the
underlying $\Q$-complexes in [De2], [Di2] without any shift of
complexes, and do not preserve mixed Hodge modules.)
The second isomorphism of (2.4.2) follows from the above argument on
the mapping cone, and the last isomorphism of (2.4.2) from
$$R\Gamma_c(X,\Q)=\R\Gamma_c(S,f_!\Q_{h,X}),\q
\R\Gamma(X,\Q)=\R\Gamma(S,f_*\Q_{h,X}).$$
\sk
Set $m=n-1$.
With the notation in the main theorem, we have the decompositions
$$\R_c\Gamma(X_s,\Q)\cong \Hsc[-m]\oplus\Q(-m)[-2m],\q
\R\Gamma(X_s,\Q)\cong \Hs[-m]\oplus\Q,$$
using the vanishing of $\Ext^i\,(i>1)$ in $\MHS$
together with the filtration $\tau_{\le k}$ as above.
We also have
$$\psi_{1/t,1}\,f_!\Q_{h,X}\cong\Hiuc[-m]\oplus\Q(-m)[-2m].$$
Let $\iota_s:\Hsc\to \Hs$ denote the canonical morphism.
The distinguished triangle (2.4.1) is then equivalent to the
isomorphism in $D^b\MHS$:
$$C\bigl(N:\Hiuc\to\Hiuc(-1)\bigr)\cong
C\bigl(\iota_s:\Hsc\to \Hs\bigr).
\leqno(2.4.3)$$
Note that $\Ker\,\iota_s$ and $\Coker\,\iota_s$ for $s\in S'$ are
extended to constant variations of mixed Hodge structures over $S$
by (1.1).
\sk
By duality we have
$$\D\bl(\Gr^W_{m-k}\Hsc\br)=
\bl(\Gr^W_{m+k}\Hs\br)(m)\q\h{for}\,\,k\ge 0.
\leqno(2.4.4)$$
Since $X_s$ is smooth affine, we have
$$\Gr^W_{m-k}\Hsc=\Gr^W_{m+k}\Hs=0
\q\h{for}\,\,k<0.$$
This implies that $\Gr^W_{m+k}\iota_s=0\,\,(k\ne 0)$, and
$\Gr^W_m\iota_s$ is an isomorphism by Proposition~2.
Combining these with the isomorphism (2.4.3) in $D^b\MHS$, we get
$$\aligned\Gr^W_{m+k}(\Ker\,N)&
\cong\begin{cases}\Gr^W_{m+k}\Hsc&\h{if}\,\,\,k<0,\\
0&\h{if}\,\,\,k\ge 0,\end{cases}\\
\Gr^W_{m+k}(\Coker\,N)&\cong
\begin{cases}\Gr^W_{m+k}\Hs&\h{if}\,\,\,k>0,\\
0&\h{if}\,\,\,k\le 0,\end{cases}\endaligned
\leqno(2.4.5)$$
where $\Coker\,N$ is a quotient of $\Hiuc(-1)$.
Using also the above duality (2.4.4), we thus get
$$\D\bl(\Gr^W_{m-k}(\Ker\,N)\br)=
\bl(\Gr^W_{m+k}(\Coker\,N)\br)(m)\q\h{for}\,\,k>0.$$
Then, applying Lemma~(2.3) to $H=\Hiuc$ where $L_kH$ is identified
with $W_k\Hsc$ for $s\in S'$, we get the second proof of Theorem~1 in
this paper.

\bs\bs
\centerline{\bf 3. Universal extensions by constant sheaves}
\bs\nin
In this section we explain the notion of a universal extension by a
constant sheaf, and then prove Theorem~3 which implies Theorem~1. 

\ms\nin
{\bf 3.1.~Universal extensions by constant sheaves over affine line.}
Let $\M$ be any pure Hodge module of weight $n$ on $S=\C$
having no constant direct factors, i.e.\ $H^{-1}(S,\M)=0$.
Note that $H^j(S,\M)=0$ for $j>0$ since $S$ is affine.
\sk
Consider the functor
$$E_{\M}(H):=\Ext^1_{\MHM(S)}(H_S[1],\M)\q\h{for}\,\,\,H\in\MHS,$$
where $H_S:=a_S^*H$ with $a_S:S\to pt$ the projection. Set
$$\HM:=H^0(S,\M),\q(\HM)_S:=a_S^*\HM.$$
Since $a_S^*$ is the left adjoint functor of $(a_S)_*$, we get the
first isomorphism of the functorial canonical isomorphisms
$$E_{\M}(H)\simto\Ext^1_{\MHS}(H[1],(a_S)_*\M)=\Hom_{\MHS}(H,\HM),
\leqno(3.1.1)$$
where the second isomorphism follows from the vanishing of $H^j(S,\M)$
for $j\ne 0$.
Note that the first isomorphism is given by taking the direct image of
$u:H_S[1]\to\M$ by $a_S$, and then composing it with the adjunction
morphism $H[1]\to(a_S)_*a_S^*H[1]$.

\sk
By (3.1.1), $E_{\M}$ is represented by $\HM$.
This means that there is a universal extension $\MM$ of $\M$ by a
constant mixed Hodge module on $S$:
$$0\to\M\to\MM\to(\HM)_S[1]\to 0,
\leqno(3.1.2)$$
whose extension class corresponds to the identity on $\HM$ by the
isomorphism (3.1.1), and such that any extension class
$\xi\in\Ext^1_{\MHM(S)}(H_S[1],\M)$ is obtained by taking the
pull-back of the short exact sequence (3.1.2) by the morphism
$$(v_{\xi})_S:H_S[1]\to(\HM)_S[1],$$
which is induced by a morphism $v_{\xi}:H\to\HM$ in $\MHS$, where
$v_{\xi}$ is uniquely determined by the extension class $\xi$.

\ms\nin
{\bf Lemma~3.2.} {\it Let
$\Hu:=\psi_{1/t,1}\M$, $\Htu:=\psi_{1/t,1}\MM$
so that we have the exact sequence in $\MHS:$
$$0\to \Hu\to\Htu\to H^{0}(S,\M)\to 0.
\leqno(3.2.1)$$
By the action of $N:=(2\pi i)^{-1}\log T_u$ together with the diagram
of the snake lemma, we have the morphism
$$\dd'':H^0(S,\M)\to\Coker(N|\Hu),
\leqno(3.2.2)$$
where $\Coker(N|\Hu)$ is a quotient of $\Hu(-1)$.
Assume the following condition holds$:$
$$\h{$\dd''$ is surjective and $\Ker\,\dd''=W_nH^0(S,\M)$.}
\leqno(C)$$
Then the weight filtration $W$ on $\Hu$ coincides with the monodromy
filtration shifted by $n$.}

\ms\nin
{\it Proof.} This follows from the primitive decomposition as in
the proof of Lemma~(1.3).

\ms\nin
{\bf 3.3.~Proof of Theorem~3.}
By Lemma~(3.2) above we have to prove condition~(C).
It is enough to show this condition for the underlying perverse
sheaves.
Let $\F$ be the underlying $\Q$-perverse sheaf of $\M$. Set
$$E_{\F}(V):=\Ext^1_{\Perv(S)}(V_S[1],\F)\q\h{for}\,\,\,
V\in M^f(\Q),$$
where $M^f(\Q)$ denotes the category of finite dimensional $\Q$-vector
spaces. By a similar argument, this functor is also represented by
$H^0(S,\F)=H^0(S,\M)_{\Q}$.
\sk
Let $\De$ be a sufficiently small open disk in $\PP^1$ with
center $\infty$ such that $\FD$ is a local system up to a
shift where $\DS:=\De\setminus\{\infty\}$.
Define for $V\in M^f(\Q)$
$$E_{\FD}(V):=
\Ext^1_{\Perv(\DS)}(V_{\DS}[1],\FD)=
\Hom_{\Q}(V,H^0(\DS,\FD)).$$
Set $\Huq:=\psi_{1/t,1}\F[-1]$.
Then $E_{\FD}$ is represented by
$$H^0(\DS,\FD)=\Coker(N:\Huq\to \Huq(-1)).$$
\sk
We have the canonical functor morphism
$$E_{\F}\to E_{\FD},$$
which corresponds to the canonical morphism
$$H^0(S,\F)=H^0(\So,\R j_*\F)\to H^0(\DS,\FD)=
H^0(\De,(\R j_*\F)|_{\De}),
\leqno(3.3.1)$$
where $j:S\into\So=\PP^1$.
We have to calculate the morphism (3.3.1).
\sk
Let $W$ be the weight filtration on $\R j_*\F$.
By the proof of [Sai2], Prop.~2.11 (see also [StZ]), we have
$$\Gr^W_k(\R j_*\F)=\begin{cases}0&\h{if}\,\,\,k<n,\\
j_{!*}\F&\h{if}\,\,\,k=n,\\
i_*\Gr^W_k\Coker(N|\Huq)&\h{if}\,\,\,k>n,\end{cases}
\leqno(3.3.2)$$
where $N:\Huq\to \Huq(-1)$ is as above, and
$i:\{\infty\}\into\So$ is the inclusion.
Since $H^{\pm 1}(\So,j_{!*}\F)=0$ by hypothesis and duality, we get
$$\Gr^W_kH^0(S,\F)=H^0(\So,\Gr^W_k(\R j_*\F)).$$

Let $j':\DS\into\De$ denote the inclusion so that
$$j_{!*}\F|_{\De}=j'_{!*}(\FD).$$
By the local classification of perverse sheaves on $\De$ (see e.g.\
[BdM], [BrMa]), we have
$$\Ext^1_{\Perv(\De)}(V_{\De}[1],j'_{!*}(\FD))=0,
\leqno(3.3.3)$$
and furthermore
$$\aligned E_{\FD}(V)
&=\Ext^1_{\Perv(\De)}(V_{\De}[1],\R j'_*(\FD))\\
&=\Ext^1_{\Perv(\De)}(V_{\De}[1],i'_*\Coker(N|\Huq))\\
&=\Hom_{\Q}(V,\Coker(N|\Huq)),\endaligned
\leqno(3.3.4)$$
where $i':\{\infty\}\into\De$.
Let $\Ft$ and $\FDt$ respectively be the universal extensions of $\F$
and $\FD$ by constant perverse sheaves so that we have the short
exact sequences in $\Perv(S)$ and $\Perv(\DS)$:
$$\aligned 0\to\F\to\Ft\to H^{0}(S,\F)_S[1]\to 0,\\
0\to\FD\to\FDt\to\Coker(N|\Huq)_{\DS}[1]&\to 0.\endaligned$$
By (3.3.3--4) we have the following commutative diagram of exact
sequences in $\Perv(\DS)$:
$$\begin{matrix}&&&&0&&0\\
&&&&\downarrow&&\downarrow\\
&&0&\to&H^{0}(\So,j_{!*}\F)_{\DS}[1]&\to&
H^{0}(\So,j_{!*}\F)_{\DS}[1]&\to&0\\
&&\downarrow&&\downarrow&&\downarrow\\
0&\to&\FD&\to&\Ft\big|_{\DS}&\to&H^{0}(S,\F)_{\DS}[1]&\to&0\\
&&||&&\downarrow&&\downarrow\\
0&\to&\FD&\to&\FDt&\to&\Coker(N|\Huq)_{\DS}[1]&\to&0\\
&&\downarrow&&\downarrow&&\downarrow\\
&&0&&0&&0&&\\ \end{matrix}\hbox{\hskip1.5cm}
\leqno(3.3.5)$$
where $\Gr^W_nH^0(S,\F)=H^0(\So,j_{!*}\F)$.
Indeed, the assertion is equivalent to that the quotient of the middle
row by the top row is isomorphic to the bottom row.
By (3.3.3--4) this follows from the fact that the restriction to
$\De$ induces the isomorphism of extension classes:
$$\Ext^1_{\Perv(\So)}(V_{\So}[1],i_*\Coker(N|\Huq))\simto
\Ext^1_{\Perv(\De)}(V_{\De}[1],i'_*\Coker(N|\Huq)).$$
\sk
Note that the morphism $\dd''$ in Lemma~(3.2) is functorially defined
for any short exact sequences on $\DS$ whose last term is constant,
and it is bijective in the case of the short exact sequence associated
to the local universal extension $\FDt$.
So the assertion follows from the above commutative diagram of
short exact sequences.

\ms\nin
{\bf 3.4.~Further property of the universal extension.}
With the notation of (3.1), let $\delta:S\into S\times S$ be the
diagonal, and $q_i:S\times S\to S$ the $i$-th projection $(i=1,2)$.
By [Sai2], 4.4.2, the inverse of the first isomorphism of (3.1.1) is
given by taking the pull-back of $v:H[1]\to(a_S)_*\M$ by $a_S$ and
then composing it with the functorial morphism:
$$a_S^*(a_S)_*\M=(q_2)_*q_1^*\M\to(q_2)_*\delta_*\delta^*q_1^*\M=\M.$$
Let $j_S$ denote the inclusion of the complement of $\delta(S)$ in
$S\times S$.
Then we have the distinguished triangle in $D^b\MHM(S)$
$$(q_2)_*(j_S)_!j_S^*q_1^*\M\to(q_2)_*q_1^*\M\to\M
\buildrel{+1}\over\to,$$
and it gives the short exact sequence (3.1.2) in $\MHM(S)$ together
with the isomorphism
$$\MM=(q_2)_*(j_S)_!j_S^*q_1^*\M[1].
\leqno(3.4.1)$$
This is essentially the same as the definition of $\MM$ by C.~Sabbah
in Appendix of [MT].
\sk
Consider the long exact sequence associated to (3.1.2):
$$0\to H^{-1}(S,\MM)\to \HM\buildrel{\dd'}\over\to
H^0(S,\M)\to H^{0}(S,\MM)\to 0,$$
where $H^{-1}(S,\M)=0$ by hypothesis and $H^j(S,\M)=H^j(S,\MM)=0$ for
$j>0$ since $S$ is affine.
Here $\dd'$ is the identity (up to a sign) by the definition the first
isomorphism in (3.1.1). So we get
$$H^j(S,\MM)=0\q\h{for any}\,\,j\in\Z.
\leqno(3.4.2)$$
\sk
Conversely, if there is a short exact sequence
$$0\to\M\to\M'\to H'_S[1]\to 0,
\leqno(3.4.3)$$
with $H'\in\MHS$ and $\M'\in\MHM(S)$ satisfying the vanishing condition
as in (3.4.2), then $\M'$ is identified with the universal extension
$\MM$ of $\M$ by a constant mixed Hodge module on $S$.
Indeed, this follows by applying the functor on the right-hand side
of (3.4.1) to the short exact sequence (3.4.3),
since $\M'=\MM'$ by the vanishing condition on $H^j(S,\M')$
and the right-hand side of (3.4.1) vanishes for a constant Hodge
module on $S$.
\sk
Since the vanishing condition as in (3.4.2) is satisfied for
$\M'=H^0f_*(\Q_{h,X}[n])$ (using the Leray spectral sequence as in
Remark~(1.2)), we get
$$\MM=H^0f_*(\Q_{h,X}[n])\q\h{with}\q\M=\Gr^W_nH^0f_*(\Q_{h,X}[n]).
\leqno(3.4.4)$$
So Theorem~3 implies the third proof of Theorem~1 in this paper.

\end{document}